\documentclass[12pt]{article}

\usepackage{cite}

\usepackage[a4paper, total={6in, 8.5in}]{geometry} 
\usepackage{hyperref} 
\usepackage{makecell} 

\usepackage{amsmath,amssymb,amsthm}
\usepackage{latexsym}
\usepackage{stmaryrd}
\usepackage{wasysym}
\usepackage{multicol}
\usepackage{multirow}
\usepackage{color}
\usepackage{dirtytalk}
\usepackage{url}
\usepackage{fontawesome} 
\usepackage{csquotes}
\usepackage{circuitikz}
\usepackage{array}

\newtheorem{theor}{Theorem}

\theoremstyle{definition}

\newtheorem{lem}[theor]{Lemma}
\newtheorem{deff}{Definition}

\newtheorem{ex}{Example}
\newtheorem{exercise}{Exercise}
\theoremstyle{remark}
\newtheorem{rem}{Remark}

\theoremstyle{definition}
\theoremstyle{definition}
\newtheorem{optimisation}{Optimisation point}

\usepackage{bm}
\def\oldvec{\mathaccent "017E\relax }
\newcommand{\Or}{{\rm O\oldvec{r}}}

\begin{document}

\pagestyle{plain} 

\title{Kontsevich graphs act on Nambu--Poisson brackets, V. Implementation}

\author{Mollie S. Jagoe Brown and  Arthemy V. Kiselev}

\maketitle

\begin{abstract}
    In this series, we established that $Q^{\gamma_3}_{d=4}(P)$ is a coboundary in 4D \cite{MSJB}, and we presented a series of experimental results about the (non)trivialisation of Kontsevich graph flows of Nambu--Poisson brackets on $\mathbb{R}^d$ \cite{IV}. This immediate sequel V. to I.--IV. \cite{avk,MSJB,fs,IV} is a guide to working with the package $\textsf{gcaops}$\footnote{\url{https://github.com/rburing/gcaops}} (\textbf{G}raph \textbf{C}omplex \textbf{A}ction \textbf{O}n \textbf{P}oisson \textbf{S}tructures) for \textsf{SageMath} by Buring \cite{BuringPhD}. Specifically, we shall explain the script used in \cite{MSJB,IV} and the use of it. 
\end{abstract}

\tableofcontents

\newpage 

\section{Introduction}

This is the fifth in the series \say{Kontsevich graphs act on Nambu--Poisson brackets} and has been written for future students who will work on the question of Kontsevich graph deformations of Poisson brackets. On its own, this paper provides all guidance needed to begin implementing the problem of establishing if $Q^{\gamma}_d(P)$ is a coboundary for Nambu--Poisson brackets $P$. The future student is encouraged to refer to the webpage of the \textsf{gcaops}\footnote{https://github.com/rburing/gcaops} software package and approach the \textsf{gcaops} developers before using it, in order to check for any important adjustments which may affect the validity of the scripts explained in this paper.

\faWarning\quad \textbf{The scripts presented here are not intended to be runnable as they are; rather they outline exactly how the scripts should be written for a given dimension $d$.}

\section{Determining if $Q^{\gamma_3}_d(P)$ is a Poisson coboundary}

We explain in detail the steps used in a script which uses \textsf{gcaops}, which can be found as an ancillary file to \cite{MSJB} on arXiv, to determine whether $Q^{\gamma}_d(P)$ is a Poisson coboundary, that is, if there exists a trivialising vector field $\vec{X}^{\gamma}_d(P)$ such that, for all Nambu--Poisson brackets $P$ on $\mathbb{R}^d$, we have that $$Q^{\gamma}_d(P)=\llbracket P,\vec{X}^\gamma_d(P)\rrbracket.$$ Our choice of cocycle $\gamma$ is the tetrahedron $\gamma_3$ used in papers I.--III. \cite{avk,MSJB,fs}, built on 4 vertices and 6 edges. In the code that will be presented, we will denote by $\textsf{p}$ the number of vertices of the cocycle $\gamma$ (\textsf{p}=\#$V(\gamma)$).

The idea is that we search for a trivialising vector field $\vec{X}^\gamma_d(P)$ over formulas with undetermined coefficients; these formulas are obtained by graphs taken with such coefficients, via the graph language due to Kontsevich, explained in Example 1 in \cite{MSJB}.

\subsection{Generate encodings}\label{generate_encodings}

It was known \cite{ascona} that the following graph gives a trivialising vector field $\vec{X}_{d=2}^{\gamma_3}(P)$:

\begin{deff}[The sunflower graph]\label{sunflower}
    The following linear combination of Kontsevich graphs (graphs built of wedges $\smash{\xleftarrow{L}\!\!\bullet\!\!\xrightarrow{R}}$, see~\cite{avk,fs,rb}) is called the sunflower graph: $$\text{sunflower }=\text{ }\raisebox{0pt}[6mm][4mm]{\unitlength=0.4mm
\special{em:linewidth 0.4pt}
\linethickness{0.4pt}
\begin{picture}(17,24)(5,5)
\put(-5,-7){
\begin{picture}(17.00,24.00)
\put(10.00,10.00){\circle*{1}}
\put(17.00,17.0){\circle*{1}}
\put(3.00,17.0){\circle*{1}}
\put(10.00,10.00){\vector(0,-1){7.30}}
\put(17.00,17.00){\vector(-1,0){14.00}}
\put(3.00,17.00){\vector(1,-1){6.67}}
\bezier{30}(3,17)(6.67,13.67)(9.67,10.33)
%
%
\put(17,17){\vector(-1,-1){6.67}}
\bezier{30}(17,17)(13.67,13.67)(10.33,10.33)
\bezier{52}(17.00,17.00)(16.33,23.33)(10.00,24.00)
\bezier{52}(10.00,24.00)(3.67,23.33)(3.00,17.00)
\put(16.8,18.2){\vector(0,-1){1}}
\put(10,17){\oval(18,18)}
\put(10,10){\line(1,0){10}}
\bezier{52}(20,10)(27,10)(21,16)
\put(21,16){\vector(-1,1){0}}
\end{picture}
}\end{picture}}=1\cdot\Gamma_1+2\cdot\Gamma_2=1\cdot \text{(0,1 ; 1,3 ; 1,2)}+2\cdot\text{(0,2 ; 1,3 ; 1,2)}.$$ The outer circle means that the outgoing arrow acts on the three vertices via the Leibniz rule. When the arrow acts on the upper two vertices, we obtain two isomorphic graphs, hence the coefficient 2 in the linear combination.
\end{deff}

Using a lucky guess in \cite{MSJB}, we took a relatively small set of graphs consisting of the sunflower graphs \ref{sunflower} which gave a trivialising vector field for the problem in $d=2$; then we \emph{expanded} \cite{MSJB,fs} the sunflower graphs to their $d=3$ and $d=4$ descendants, and we searched for solutions $\vec{X}^{\gamma_3}_{d=3}(P)$ and $\vec{X}^{\gamma_3}_{d=4}(P)$ over the formulas obtained from such graphs, with undetermined coefficients. In this way, we discovered in \cite{MSJB} that over dimension $d=4$, there exists a trivialising vector field $\vec{X}^{\gamma_3}_{d=4}(P)$ such that $Q^{\gamma_3}_{d=4}(P)=\llbracket P,\vec{X}^{\gamma_3}_{d=4}(P)\rrbracket$. 

\begin{ex}
    Let us show how one obtains the encodings of the $d=3$ and $d=4$ descendants from the encodings of the sunflower graph (Definition \ref{sunflower}) which gives a trivialising vector field $\vec{X}^{\gamma_3}_d(P)$ for the problem in $d=2$ (see section 2.2, Definition 6 and Example 4 in \cite{MSJB}) for details. Essentially, we plug the dimension-specific Nambu--Poisson structures into the sunflower's vertices. We give the graphical intuition in Table \ref{buildingblocks}, showing the graphical expression of the dimension-specific Nambu--Poisson structures with which we build the sunflower graphs in dimensions $d=3$ and $d=4$. 
    
    \begin{table}\caption{Number of sunflower graphs in dimensions $d=2,3,4$ due to the Nambu--Poisson building blocks.}\label{buildingblocks}
\begin{tabular}{l|l|l|l|}
\cline{2-4}
 & 2D & 3D & 4D \\ \hline
\multicolumn{1}{|l|}{Building blocks} & \begin{minipage}{0.2\textwidth}\includegraphics[width=0.8\linewidth]{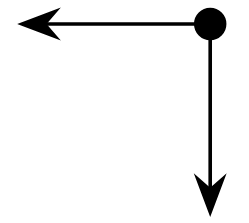}\end{minipage} & \begin{minipage}{0.2\textwidth}\includegraphics[width=\linewidth]{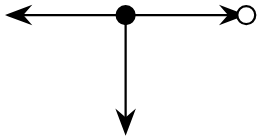}\end{minipage} & 
\begin{minipage}{0.2\textwidth}\includegraphics[width=\linewidth]{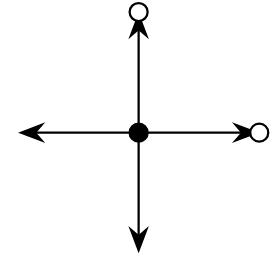}\end{minipage}
 \\ \hline
\multicolumn{1}{|l|}{Number of sunflower graphs} & 3 & 48 & 324 \\ \hline
\end{tabular}
\end{table}

    In the sought 1-vectors $\vec{X}_d^{\gamma_3}(P)$ each graph is built of three Nambu--Poisson structures over $\mathbb{R}^d$. Consider the $d=2$ encoding of the first graph of the sunflower:
    $$ (0,1;1,3;1,2).$$ This encoding can be read as follows. The semi-colon separates each Nambu--Poisson structure $\smash{\xleftarrow{L}\!\!\bullet\!\!\xrightarrow{R}}$ indicating which arrow (left or right) goes where. The Nambu--Poisson structure based at vertex 1 has its left arrow acting on vertex 0 (the sink) and its right arrow acting on itself. The Nambu--Poisson structure based at vertex 2 has its left arrow acting on vertex 1 and its right arrow acting on vertex 3. The Nambu--Poisson structure based at vertex 3 has its left arrow acting on vertex 1 and its right arrow acting on vertex 2.

    To obtain the encodings of this $d=2$ graph's 3D-descendants, we insert the $d=3$ Nambu--Poisson bi-vectors in each vertex, and expand via the Leibniz rule. The encodings are given by $$(0,i_1,4;i_2,j,5;i_3,k,6)$$ for $i_1,i_2,i_3\in\{1,4\},j\in\{3,6\},k\in\{2,5\}$. To generate all encodings of 3D-descendants, we use \textsf{SageMath} code: \begin{verbatim}
X_graph_encodings = []
for i2 in [2,5]:
   for j1 in [1,4]:
       for j2 in [1,4]:
           for k in [3,6]:
               X_graph_encodings.append((0,1,4,j1,k,5,j2,i2,6))
    \end{verbatim}
\end{ex} \hfill $\blacksquare$

\noindent Recall that the tuples separated by a semi-colon each pertain to a copy of the Nambu--Poisson bracket; and this number is distinct for the problem $Q^\gamma_d(P)=\llbracket P,\vec{X}^\gamma_d(P)\rrbracket$ for each $\gamma$. Note that we are not able to use semi-colons in the code to separate the tuples given by the Nambu--Poisson structures. We will truncate the encodings accordingly in the next step.

\begin{exercise}
    Show why the code to generate the encodings of the 4D-descendants of the first graph of the sunflower is given by:
    \begin{verbatim}
X_graph_encodings = []
for i2 in [2,5,8]:
   for j1 in [1,4,7]:
       for j2 in [1,4,7]:
           for k in [3,6,9]:
               X_graph_encodings.append((0,1,4,7,j1,k,5,8,j2,i2,6,9))
    \end{verbatim}

    \noindent \textit{Hint:} The new Casimirs appearing in the 4D Nambu--Poisson bracket can be numbered by 7,8,9, in order. That is, vertex 1 has Casimirs $a_1$ denoted by 4, and $a_2$ denoted by 7. Vertex 2 has Casimirs $a_1$ denoted by 5, and $a_2$ denoted by 8. Vertex 3 has Casimirs $a_1$ denoted by 6, and $a_2$ denoted by 9.  
\end{exercise} 

\subsection{Convert Kontsevich graph encodings to Nambu micro-graphs}\label{encodings_to_mgs}

We now construct the graphs from the encodings we generated in \S\ref{generate_encodings}. We define the mapping \verb|encoding_to_graph| which does this operation. To begin, we truncate the encodings of the $d$-descendant graphs into tuples, where each tuple corresponds to a copy of the Nambu--Poisson structure. In the code, define \verb|target| in such a way that \verb|targets[i]| is the list of vertices that the $i$-th Levi-Civita vertex points to; that is, \verb|targets[i]| refers to the targets of the $i$-th Levi-Civita vertex. As we saw in \S\ref{generate_encodings} for $X^{\gamma_3}$, we are using encodings of graphs built from 3 Nambu--Poisson structures, hence we will divide each encoding into 3 tuples. Generally, we divide $d$-dimensional graph encodings into as many tuples as there are vertices in the graphs over which we search for a trivialising vector field. \begin{verbatim}
    from gcaops.graph.formality_graph import FormalityGraph
    def encoding_to_graph(encoding):
    targets = [encoding[0:d], encoding[d:2d], encoding[2d:3d]]
\end{verbatim}

\noindent We now have our truncated encodings, and we define the edges of the graphs we will assemble. \begin{verbatim}
edges = sum([[(k+1,v) for v in t] for (k,t) in enumerate(targets)], [])
\end{verbatim}

\noindent We have that \verb|enumerate(targets)| gives elements of the form \verb|(k,t)| where $k$ enumerates the tuples $t$ in \verb|targets|. Here, we have that \verb|enumerate(targets)| gives $$\bigl(0,t=(0,1,4)\bigr)\quad  \bigl(1,t=(1,3,5)\bigr)\quad \bigl(2,t=(1,2,6)\bigr).$$ The edges can then be created using \begin{verbatim}     [[(k+1,v) for v in t] for (k,t) in enumerate(targets)]\end{verbatim} Here, for $k=0$ this looks like: $$(1,0)\quad (1,1) \quad (1,4);$$ for $k=1$ this looks like: $$(2,1)\quad (2,3)\quad(2,5);$$ for $k=2$ this looks like $$(3,1)\quad(3,2)\quad(3,6).$$ The edges are then collected in a list with \verb|sum|. 

Finally, we assemble the list of Nambu--Poisson \verb|Formality graphs|. These are given in terms of $m$ sinks (a function can be plugged into each sink -- in this way, each \verb|Formality graph| with $m$ sinks represents a $m$-vector, by applying Kontsevich's graph language, see example 1 in \cite{MSJB}), $n$ aerial vertices, and a set of edges. \begin{verbatim}
    return FormalityGraph(m, n, edges)
\end{verbatim}

\begin{ex}
    In our 2D and 3D, we have $$\text{2D}: \verb|FormalityGraph(1,3,edges)| \quad \text{3D}: \verb|FormalityGraph(1,6,edges)|$$ Note that for the 2D graphs which give a trivialising vector field, we have $m=1$ in any dimension, because its graphs represent a vector field, and $n=3+3(d-2)$, where $d-2$ is the number of Casimirs in the Nambu--Poisson structure over $\mathbb{R}^d$. 
\end{ex} \hfill $\blacksquare$

Overall, we now have the function \verb|encoding_to_graph|: \begin{verbatim}
def encoding_to_graph(encoding):
    targets = [encoding[0:d], encoding[d:2d], encoding[2d:3d]]
    edges = sum([[(k+1,v) for v in t] for (k,t) in enumerate(targets)], [])
    return FormalityGraph(m, n, edges)
\end{verbatim} We can generate the Nambu--Poisson Formality graphs that we need: \begin{verbatim}
    X_graphs = [encoding_to_graph(e) for e in X_graph_encodings]
    print("Number of graphs in X_graphs matches X_graph_encodings:", 
            len(X_graphs)==len(X_graph_encodings), flush=True)
\end{verbatim}

\begin{rem}
    Since the results from \cite{MSJB} were obtained, a built-in isomorphism checker has been added to the package \textsf{gcaops}. Isomorphic graphs give the same formula, so this increases efficiency due to reducing the number of graphs we use from this step onwards.
\end{rem}

\subsection{How to obtain formulas from micro-graphs}\label{graph_to_fla_eval}

We work with formulas obtained from the graphs produced in \S\ref{encodings_to_mgs}, by using Kontsevich's graph language. The key premise is that the edges are directed arrows which contain derivatives $$(\text{tail})\frac{\overleftarrow{\partial}}{\partial\xi_i}\otimes\frac{\overrightarrow{\partial}}{\partial x^i}(\text{head})$$ so when an arrow hits a vertex, it in fact acts as a derivative on the target vertex content. An example is given in Example 1, \cite{MSJB}.

We begin by setting up the necessary differential polynomial ring $D_d$ and superfunction algebra $S_d$, where $d$ is the base dimension. \begin{verbatim}
    from gcaops.algebra.differential_polynomial_ring import 
            DifferentialPolynomialRing
    Dd = DifferentialPolynomialRing(QQ, ('rho','a1','a2',...,'a(d-2)'),
                ('x0','x1',...,'x(d-1)'), 
            max_differential_orders=[d,d+1,d+1])
    rho, a1, a2,..., a(d-2) = Dd.fibre_variables()
    x0,x1,...,x(d-1) = Dd.base_variables()
    even_coords = [x0,x1,...,x(d-1)]

    from gcaops.algebra.superfunction_algebra import SuperfunctionAlgebra
    Sd.<xi0,xi1,...,xi(d-1)> = SuperfunctionAlgebra(Dd, Dd.base_variables())
    xi = Sd.gens()
    odd_coords = xi
    epsilon = xi[0]*xi[1]*...*xi[d-1] # Levi-Civita tensor
\end{verbatim}

\noindent The line \begin{verbatim}
Sd.<xi0,xi1,...,xi(d-1)> = SuperfunctionAlgebra(Dd, Dd.base_variables())
\end{verbatim} means that parity-odd coordinates $\xi_{i}$ correspond to the parity-even coordinates $x_i$. In dimension $d=4$ we would have for instance $(x^0,x^1,x^2,x^3)=(x,y,z,w)$. The Levi-Civita symbol $\varepsilon^{i_1i_2\hdots i_d}$ is \verb|epsilon|; it helps us to encode the determinant in the Nambu--Poisson bracket. 

Let us create the function that takes Nambu--Poisson Formality graphs and gives their $k$-vector formulas. 

\begin{rem}
    The package \textsf{gcaops} uses built-in multi-processing. To use multi-processing again in the code would result in an error. 
\end{rem}

We use a trick: create the Euler vector field $E=\sum_{i=0}^{d-1}x^i\xi_{i}$ and insert it into the sink. In this way, the incoming derivation $\partial/\partial x_i$ results in $\xi_i$, for $i$ running from 1 to $d$.

\begin{verbatim}
    E = x0*xi[0] + x1*xi[1] + ... + x(d-1)*xi[d-1]
\end{verbatim}

The function initially sets its future output equal to zero, in the superfunction algebra $S_d$. It proceeds to add the computed formulas of graphs to it.

\begin{verbatim}
    result = Sd.zero()
\end{verbatim}

\noindent We begin a for-loop over the indices labelling the edges of the micro-graphs. 

\begin{verbatim}
for index_choice in itertools.product(itertools.permutations(range(d)),
            repeat=p):
\end{verbatim}

\noindent Here, \verb|itertools.permutations(range(d))| is just $(S_{\{1,2,...,d\}})^p$, so each element is a $p$-tuple of permutations of $\{1,...,d\}$. Each element \verb|index_choice| looks like $[\sigma_1(i),\sigma_2(i),\sigma_3(i)]$ with $\sigma_1,\sigma_2,\sigma_3\in S_{\{1,2,...,d\}}$ and $i\in\{1,...,d\}$. So, each $\sigma_j(i)$ is a string of $d$ indices, for instance \verb|('1','2',...,'d')|. Note that we have $p=3$ elements $\sigma_j(i)$ in \verb|index_choice| because the micro-graphs we use for $\vec{X}_d^{\gamma_3}(P)$ are built of $p=3$ Nambu--Poisson structures. Namely, we want to obtain all $d!$ possible permutations of outgoing edges in each Nambu--Poisson vertex. For micro-graphs in dimension $d$ built of $p$-many Nambu--Poisson brackets, we have $(d!)^p$ many \verb|index_choice| elements.


\textit{Hint:} This can be used to check if your code is correct. For instance, the script could contain a line giving the condition that if the length of the object \newline \verb|itertools.product(itertools.permutations(range(d)), repeat=p)| is not equal to $(d!)^p$, that is, \begin{verbatim}len(itertools.product(itertools.permutations(range(d)), repeat=p)
        !=factorial(d)**p,\end{verbatim} where \verb|!=| means $\neq$ and \verb|factorial(d)**p| means $(d!)^p$, then the code stops running. 

\begin{ex}
    We take the example of the 3D micro-graph $\Gamma$ given by the encoding $(0,2,4;1,3,5;1,2,6)$. We write the inert sum of its formula, and find that we can rephrase the summation over indices by summation over permutations $\sigma_i\in S_{\{1,2,3\}}$: 
    \begin{align*}
        \phi(\Gamma) &=\sum_{\vec{\imath},\vec{\jmath},\vec{k}=1}^{d=3}\varepsilon^{i_1i_2i_3}\varepsilon^{j_1j_2j_3}\varepsilon^{k_1k_2k_3}\varrho_{j_1k_1}\varrho_{i_2k_2}\varrho_{j_2}a_{i_3}a_{j_3}a_{k_3}\partial_{i_1}() \\
        &= \sum_{\sigma_1,\sigma_2,\sigma_3\in S_{\{1,2,3\}}} (-)^{\sigma_1}(-)^{\sigma_2}(-)^{\sigma_3}\varrho_{\sigma_2(1)\sigma_3(1)}\varrho_{\sigma_1(2)\sigma_3(2)}\varrho_{\sigma_2(2)}a_{\sigma_1(3)}a_{\sigma_2(3)}a_{\sigma_3(3)}\partial_{\sigma_1(1)}().
    \end{align*}
    The indices $\vec{\imath}=i_1,i_2,i_3$ corresponds to $\sigma_1$, so $\sigma_1(\alpha)$ corresponds to the index $i_\alpha$. Similarly, the indices $\vec{\jmath}=j_1,j_2,j_3$ corresponds to $\sigma_2$ and $\vec{k}=k_1,k_2,k_3$ corresponds to $\sigma_3$. 
\end{ex} \hfill $\blacksquare$ 

\noindent Recalling that $\varepsilon^{i_1i_2i_3}\in\{-1,0,1\}$, we know that $\varepsilon^{\sigma(i_1i_2i_3)}\varepsilon^{\sigma(j_1i_2j_3)}\varepsilon^{\sigma(k_1k_2k_3)}\in\{-1,0,1\}$ for any $\sigma\in S_{\{1,2,...,d\}}$. We create the \verb|sign| in the final sum:
\begin{verbatim}sign=epsilon[index_choice[0]]*epsilon[index_choice[1]]*epsilon[index_choice[2]]
\end{verbatim}

If our micro-graphs were built of $p$-many Nambu--Poisson structures, then the last term above would end up with \verb|epsilon[index_choice[p-1]]|.

The following assumes that the ground vertex is labelled as 0; the aerial vertices of out-degree $d$ are labelled 1,2,...,$d$; the $d$-many Casimirs $a^1$ are labelled $d+1,...,2d$;...; the $d$-many Casimirs $a^{d-2}$ are labelled $(d-2)d+1,...,(d-1)d$.  

\begin{ex}
    In 3D, the ground vertex is 0; the aerial vertices of out-degree 3 are 1,2,3; the Casimirs $a^1$ are 4,5,6. 

    In 4D, we have the same as above in 3D, with the added Casimirs $a^2$ being 7,8,9.
\end{ex}\hfill$\blacksquare$

We define the vertex content of the micro-graphs. Recall that the micro-graphs we use are built of three (this is specific for the problem with $\gamma_3$) Nambu--Poisson bi-vectors $P(\boldsymbol{a},\varrho)$. That is, the vertex content can be understood as containing the smooth functions $\varrho,a^1,a^2,...,a^{d-2}$. 

\begin{verbatim}
vertex_content = [E, Sd(rho), Sd(rho), Sd(rho), Sd(a1), Sd(a1), 
        Sd(a1), Sd(a2), Sd(a2), Sd(a2),...,Sd(a(d-2)),Sd(a(d-2)),Sd(a(d-2))]
\end{verbatim}

We take the set \verb|g.edges()| of edges of the graph. They are of the form \verb|(source,target)|, meaning that the edge is directed from \verb|source| to \verb|target|. We assign an index to each edge using \verb|zip|: 
\begin{verbatim}
    for ((source, target), index) in zip(g.edges(), 
        sum(map(list, index_choice), [])):  
\end{verbatim} where \verb|map(list, index_choice)| creates the elements of \verb|index_choice| as separate lists; and \verb|sum(map(list,index_choice), []))| adds all \verb|index_choice| elements as one big list. 

\begin{ex}
Let us take a toy example to visualise the meaning from the above:
\begin{verbatim}
    index_choice=[('1','2','3'),('4','5','6')]
\end{verbatim}
Applying \verb|list| to this, we have:
\begin{verbatim}    print(list(index_choice)): [('1', '2', '3'), ('4', '5', '6')]\end{verbatim}
But taking \verb|index_choice| as a \verb|list| under the built-in \verb|map| just gives the dictionary (hash map): \begin{verbatim}    print(map(list,index_choice)): <map object at 0x14af1bf5f160> \end{verbatim}
But if we \verb|sum| this to an empty list \verb|[]|, we obtain what we want: \begin{verbatim}print(sum(map(list, index_choice), [])): ['1', '2', '3', '4', '5', '6']\end{verbatim}
This is equivalent to applying \verb|list| once more:
\begin{verbatim}print(list(map(list, index_choice))): ['1', '2', '3', '4', '5', '6']\end{verbatim}
\end{ex} \hfill $\blacksquare$

\begin{rem}
    The takeaway is that the built-in function \verb|map| does not immediately give output on its own, one needs to pick what type of object should be output. In the current code, as in the example above, the chosen object to be output is a list. 
\end{rem}

The vertex content is differentiated, which is due to the action of the incoming arrows: 
\begin{verbatim}
    vertex_content[target] = vertex_content[target].derivative(
            even_coords[index]) 
\end{verbatim}

The resulting output of the function requires the final assembly of the sum, meaning multiplication of the content of the vertices by \verb|sign|:

\begin{verbatim}
   result += sign * prod(vertex_content) 
\end{verbatim}

Mathematically, this concludes the process going from a graph to its superalgebra representation of a 1-vector formula. 

Assuming that the ground vertex is labeled 0, the vertices of out-degree 4 are labeled 1, 2, 3, and the Casimirs are labeled 4, 5, 6; 7, 8, 9, the full function which takes 4D graphs as input and outputs their 1-vector formulas is given by: 

\begin{verbatim}
def evaluate_graph(g):
    E = x*xi[0] + y*xi[1] + z*xi[2] + w*xi[3] 
    # Euler vector field, to insert into ground vertex. 
    #Incoming derivative d/dx^i will result in xi[i].
    result = S4.zero()
    for index_choice in itertools.product(itertools.permutations(range(
                    4)), repeat=3):
        sign = epsilon[index_choice[0]] * epsilon[index_choice[1]] * 
                epsilon[index_choice[2]]
        vertex_content = [E, S4(rho), S4(rho), S4(rho), S4(a1), 
                        S4(a1), S4(a1), S4(a2), S4(a2), S4(a2)]
        for ((source, target), index) in zip(g.edges(), 
                    sum(map(list, index_choice), [])):
            vertex_content[target] = vertex_content[target].derivative(
                        even_coords[index])
        result += sign * prod(vertex_content)
    return result

\end{verbatim}

\begin{rem}
    The micro-graphs we work with give 1-vector formulas, but this code can be adapted to give other $k$-vector formulas from suitable micro-graphs. 
\end{rem}

We can now use the function defined above to compute the formulas from the micro-graphs. The built-in function \verb|imap| is useful for this step because it returns jobs as soon as they are done, instead of waiting for all jobs to be done (which is the case of \verb|map|). Thus, we have the following code to obtain the formulas:

\begin{verbatim}
    X_formulas = []
    with Pool(processes=NPROCS) as pool:
        X_formulas = list(pool.imap(evaluate_graph, X_graphs))
        X_formulas = [S4(X_formula) for X_formula in X_formulas]
\end{verbatim}

Above, we apply the function \verb|evaluate_graph| to all elements of the list \verb|X_graphs| using \verb|imap|; this is executed with \verb|Pool|. The resulting formulas are stored as a list \verb|X_formulas|. In the last line, the call \verb|Sd(X_formula)| converts the plain-text formulas into usable elements of the superfunction algebra $S_d$. 

\begin{optimisation} This method consumes large amounts of memory; this could be optimised. We believe that the following method would be efficient in reducing the amount of memory used. Compute the formulas one by one, and sum them together with undetermined coefficients instantly, so that similar terms are collected. We then have formulas of the form $(c_1+c_2)\cdot\text{term}$ instead of $c_1\cdot\text{term}+c_2\cdot\text{term}$. 
\end{optimisation}

\subsection{Collect monomials from formulas}\label{collect_monomials}

We now aim to create a basis of monomials in the formulas obtained in \S\ref{graph_to_fla_eval}. Inputting a micro-graph $g$ into the function described in \S\ref{graph_to_fla_eval} gives a formula $\phi(g)$ of the form $$\phi(g)=\sum_{i=0}^{d-1} f_g^i\xi_i,$$ where $f_\Gamma^i$ is the differential polynomial coefficient of $\xi_i$ depending on the graph $g$. We want to collect the distinct differential monomials in the formula. Namely, for each $i=1,2,...,d$ we collect the distinct differential monomials in $f_g^i$. Specifically, we obtain a basis of such monomials. We begin by giving a simple example of this process in 3D. 

\begin{ex}
Let us assume that a micro-graph $g_1$ gives the following 3D formula: $$\phi(g_1)=\bigl(\varrho_xa_{xx}+3a_z\bigr)\xi_0 + \bigl(\varrho_{zy}a_{y}-a_z+2\varrho_{xyz}a_{x}\bigr)\xi_2.$$ For $i=1,2,3$ we can identify its components $f_{g_1}^i$ as follows: $$f^1_{g_1}=\varrho_xa_{xx}+3a_z,\quad f^2_{g_1}=0, \quad f^3_{g_1}=\varrho_{zy}a_{y}-a_z+2\varrho_{xyz}a_{x}.$$ Then the distinct non-zero differential monomials near each $\xi_i$ in the formula $\phi(g_1)$ are: 
$$\xi_0:\varrho_xa_{xx},3a_z,\quad \xi_1:\varnothing,\quad \xi_2:\varrho_{zy}a_{y},-a_z,2\varrho_{xyz}a_{x}.$$ Naturally, bases of the above distinct differential monomials near each $\xi_i$ is given by: $$\xi_0:\varrho_xa_{xx},a_z,\quad \xi_1:\varnothing,\quad\xi_2:\varrho_{zy}a_{y},a_z,\varrho_{xyz}a_{x}.$$ We see that for $\xi_0$, the differential monomial $3a_z$ can indeed be obtained from the basis element $a_z$, and for $\xi_2$, the differential monomial $-a_z$ can each be obtained by the basis element $a_z$.
\end{ex}\hfill$\blacksquare$

To execute this procedure in \textsf{SageMath}, we first collect the distinct differential monomials for each $\xi_i$. To this end, we define a list of $d$-many sets. Each $i$-th set will contain all distinct differential monomials of all formulas in the coefficient of $\xi_i$. We use \verb|set| because it does not count duplicates, so we can indeed collect distinct monomials. 

\begin{verbatim}
   X_monomial_basis = [set([]) for i in range(d)] 
\end{verbatim}

We now create nested for-loops in order to pick the differential polynomial coefficient $f_g^i$ of each $\xi_i$ for each formula $\phi(g)$. The first for-loop will select $i\in\{0,1,2,...,d-1\}$ to pick a certain $\xi_{i}$, the second for-loop will select a formula in the list of formulas obtained from micro-graphs \verb|X_formulas|. 

\begin{verbatim}
    for i in range(d):
        for X in X_formulas:
            X_monomial_basis[i] |= set(X[i].monomials())
\end{verbatim}

Here, we collect the monomials in each \verb|X[i]| term of the formula \verb|X|, which corresponds to collecting the monomials in each $f_g^i$ term of the formula $\phi(g)$. This is done using the built-in command \verb|.monomials()|. So, \verb|X[i].monomials()| identifies all monomials in \verb|X[i]|; then because \verb|set()| does not allow duplicates, \verb|set(X[i].monomials())| gives the distinct monomials in \verb|X[i]|. We add the distinct monomials in \verb|X[i]| to \verb|X_monomial_basis[i]|, which is the $i$-th set in \verb|X_monomial_basis|. 

\begin{rem}
    We emphasise that we collect the distinct monomials \emph{for each} $\xi_i$. That is, we could have the same monomial appear in \verb|X_monomial[0]|, which corresponds to $\xi_0$, and \verb|X_monomial[3]|, which corresponds to $\xi_4$. 
\end{rem}

\begin{rem}
    The command \verb|.monomials()| only recognises non-zero terms as monomials. Therefore, if we have for instance \verb|X[0]| equal to 0, then \verb|X[0].monomials()| returns \verb|()|; that is, nothing is done. If we have \verb|X[0]| equal to 1, then \verb|X[0].monomials()| returns \verb|(1,)|.
\end{rem}

\begin{rem}
    \faWarning The command \verb|.monomials()| does \emph{not} take into account coefficients. For example, if we have $f=2xy+4y$, then \verb|f.monomials()| returns $[xy,y]$. The coefficients can be obtained by the command \verb|.coefficients()|, so here we have that \verb|f.coefficients()| returns $[2,4]$. This means that for us, the way \verb|X_monomial_basis| is defined indeed provides us with a true basis of differential monomials for each $\xi_{i_1}$. 
\end{rem}

We then store the sets in \verb|X_monomial_basis| as lists. This will allow us to enumerate the monomials in the next step.

\begin{verbatim}
    X_monomial_basis = [list(b) for b in X_monomial_basis]
\end{verbatim}

\noindent We assign an index to each monomial:

\begin{verbatim}
    X_monomial_index = [{m : k for k, m in enumerate(b)} for b 
                            in X_monomial_basis] 
\end{verbatim}
The built-in command \verb|enumerate()| assigns an index to each element of the lists in \verb|X_monomial_basis|. However, this enumeration cannot be displayed immediately, so that the monomial $m$ with index $k$ is added to a dictionary (hash map) with key $m$ and value $k$. This is written as $\{m:k\}$. Note that \verb|X_monomial_index| will be a list of $d$-many collections of indexed monomials. The point of this step is to be able to obtain the index of a differential monomial in a very fast way, instead of iterating through the whole list of differential monomials each time.

Let us exemplify all of the steps above using three arbitrary formulas in 3D.
\begin{ex}
We take three 3D formulas assumed to be obtained from three distinct non-isomorphic micro-graphs $g_1,g_2,g_3$:
    \begin{align*}
\phi(g_1)&=\bigl(\varrho_xa_{xx}+3a_z\bigr)\xi_0 + \bigl(\varrho_{zy}a_{y}-a_z+2\varrho_{xyz}a_{x}\bigr)\xi_2,\\
\phi(g_2)&=\bigl(a_z\bigr)\xi_0 + \bigl(\varrho_xa_{xx}\bigr)\xi_1+\bigl(2\varrho_{zy}a_{y}+\varrho_{xyz}a_{x}\bigr)\xi_2, \\
\phi(g_3)&=\bigl(\varrho_xa_{xx}+a_z\bigr)\xi_0+\bigl(-2\varrho_{x}a_{xx}\bigr)\xi_1 + \bigl(-3\varrho_{zy}a_{y}-a_z\bigr)\xi_2.
\end{align*}

\noindent After the nested for-loops, we have \verb|X_monomial_basis| given by: $$[\{\varrho_xa_{xx},a_z\}\quad \{ \varrho_{x}a_{xx}\} \quad \{\varrho_{zy}a_{y},a_z,\varrho_{xyz}a_x\}].$$ When storing the sets in \verb|X_monomial_basis| as lists, we obtain: $$[[\varrho_xa_{xx},a_z]\quad [ \varrho_{x}a_{xx}] \quad [\varrho_{zy}a_{y},a_z,\varrho_{xyz}a_x]].$$ 

\noindent Taking \verb|b in X_monomial_basis| means taking one of the three lists in \verb|X_monomial_basis|. Recall that each list corresponds to $\xi_{i_1}$. Then, the built-in command \verb|enumerate()| assigns an index to each element of the lists. We choose to display the enumeration in the form $\{m:k\}$, where $m$ corresponds to the monomial, and $k$ corresponds to the index. In our example, using the first list in \verb|X_monomial_basis|, we have: 
$$b=[\varrho_xa_{xx},a_z],$$ then \verb|{m : k for k, m in enumerate(b)}| gives $$\{[\varrho_xa_{xx}:0,\quad a_z:1]\}.$$ Therefore, \verb|X_monomial_index| is given by: $$[\{\varrho_xa_{xx}:0,\quad a_z:1\}\quad \{ \varrho_{x}a_{xx}:0\} \quad \{\varrho_{zy}a_{y}:0,\quad a_z:1,\quad\varrho_{xyz}a_x:2\}].$$ 

\noindent Note that the count of indices starts again for each $\xi_i$ because we are collecting the distinct differential monomials \emph{for each} $\xi_i$. We store the total number of monomials as \verb|X_monomial_count|:
\begin{verbatim}
   X_monomial_count = sum(len(b) for b in X_monomial_basis) 
\end{verbatim}
In this example, we have \verb|X_monomial_count|=6. 
\end{ex}\hfill$\blacksquare$

Overall, the code for collecting the distinct differential monomials contained in all formulas obtained in \S\ref{graph_to_fla_eval} is as follows. 

\begin{verbatim}
X_monomial_basis = [set([]) for i in range(d)]
for i in range(d):
    for X in X_formulas:
        X_monomial_basis[i] |= set(X[i].monomials())
X_monomial_basis = [list(b) for b in X_monomial_basis]
X_monomial_index = [{m : k for k, m in enumerate(b)} for b in X_monomial_basis]

print("Number of monomials in components of X:", 
        [len(b) for b in X_monomial_basis], flush=True)

X_monomial_count = sum(len(b) for b in X_monomial_basis)
\end{verbatim}
In the next step, we will use our new basis of monomials in order to express the graph-to-formula evaluation as a matrix. 

\subsection{Express graph-to-formula evaluation as a matrix}\label{eval_matrix}

So far, we have converted encodings to micro-graphs in \S\ref{encodings_to_mgs}, converted micro-graphs to formulas in \S\ref{graph_to_fla_eval} and collected the monomials near each $\xi_i$, $i=0,1,2,...,d-1$ in the formulas in \S\ref{collect_monomials}. We will now express the graph-to-formula evaluation as a matrix; that is, we will obtain a linear map from coefficients of graphs to the coefficients of the differential monomials in a suitable basis. This will allow us to identify which graphs are linearly independent as formulas. We begin by showing this matrix for the example in \S\ref{collect_monomials}. The authors highly encourage the reader to rely on this example for intuition and understanding of this step of the code. 

\begin{ex}
   We take three 3D formulas assumed to be obtained from three distinct non-isomorphic micro-graphs $g_1,g_2,g_3$. 
    \begin{align*}
\phi(g_1)&=\bigl(\varrho_xa_{xx}+3a_z\bigr)\xi_0 + \bigl(\varrho_{zy}a_{y}-a_z+2\varrho_{xyz}a_{x}\bigr)\xi_2,\\
\phi(g_2)&=\bigl(a_z\bigr)\xi_0 + \bigl(\varrho_xa_{xx}\bigr)\xi_1+\bigl(2\varrho_{zy}a_{y}+\varrho_{xyz}a_{x}\bigr)\xi_2, \\
\phi(g_3)&=\bigl(\varrho_xa_{xx}+a_z\bigr)\xi_0+\bigl(-2\varrho_{x}a_{xx}\bigr)\xi_1 + \bigl(-3\varrho_{zy}a_{y}-a_z\bigr)\xi_2.
\end{align*} We saw in \S\ref{collect_monomials} that the corresponding basis of this set of formulas is given by: $$[[\varrho_xa_{xx},a_z]\quad [ \varrho_{x}a_{xx}] \quad [\varrho_{zy}a_{y},a_z,\varrho_{xyz}a_x]].$$ The matrix of coefficients that we construct can be understood as follows, with columns corresponding to each $\phi(g_k)$: $$\begin{bmatrix}
    \varrho_xa_{xx}\text{ near }\xi_0 \\
    a_z\text{ near }\xi_0 \\
    \varrho_xa_{xx}\text{ near }\xi_1\\
    \varrho_{zy}a_y\text{ near }\xi_2 \\
    a_z\text{ near }\xi_2 \\
    \varrho_{xyz}a_x\text{ near }\xi_2
\end{bmatrix}.$$ In our example, we have three formulas from three micro-graphs $g_1,g_2,g_3$. Each formula corresponds to a column in the matrix. Therefore, our matrix is: $$\begin{bmatrix}
    1&0&1\\3&1&1\\0&1&-2\\1&2&-3\\-1&0&-1\\2&1&0
\end{bmatrix}.$$ The first column can be read as follows: in $\phi(g_1)$, the coefficient of $\varrho_xa_{xx}$ near $\xi_0$ is 1; the coefficient of $a_z$ near $\xi_0$ is 3; the coefficient of $\varrho_xa_{xx}$ near $\xi_1$ is 0; the coefficient of $\varrho_{zy}a_y$ near $\xi_2$ is 1; the coefficient of $a_z$ near $\xi_2$ is $-1$; the coefficient of $\varrho_{xyz}a_x$ near $\xi_2$ is 2. 

\begin{exercise}
    Check the other two columns.
\end{exercise}

We have obtained the graph-to formula evaluation matrix of our example with three micro-graphs $g_1,g_2,g_3$. 
\end{ex}\hfill$\blacksquare$

Now, we turn back to the explanation of the code for obtaining the graph-to-formula evaluation matrix of the general problem. Essentially, we define a matrix which has as many columns as we have formulas, and as many rows as we have monomials \emph{for each} $\xi_i$. The number of rows is given by \verb|X_monomial_count|. We then fill the matrix entries with the monomial coefficients, \emph{looking at each $\xi_i$ term for each formula}. We create the individual columns, each column corresponding to a formula.  

\begin{verbatim}
X_monomial_count = sum(len(b) for b in X_monomial_basis)

X_evaluation_matrix = matrix(QQ, X_monomial_count, 
                        len(X_formulas), sparse=True)
\end{verbatim}
We now begin filling the matrix entries. To do this, we create empty vectors, which will become the columns of the matrix. We need two for-loops: the first will pick a formula from \verb|X_formulas|, the second will pick a $\xi_i$ term of the formula.
\begin{verbatim}
for i in range(len(X_formulas)):
    v = vector(QQ, X_monomial_count, sparse=True)
    index_shift = 0
    for j in range(d):
        f = X_formulas[i][j]
        for coeff, monomial in zip(f.coefficients(), f.monomials()):
            monomial_index = X_monomial_index[j][monomial]
            v[index_shift + monomial_index] = coeff
        index_shift += len(X_monomial_basis[j])
    X_evaluation_matrix.set_column(i, v)
\end{verbatim}
Recall that each column is truncated into \say{boxes} by the distinct $\xi_i$ terms. To account for this when filling the columns, we use \verb|index_shift| in order to ensure that the coefficients of monomials go to the correct $(\xi_i)$--boxes. Indeed, the second for-loop \verb|for j in range(d)| allows us to pick \verb|X_formulas[i][j]|, that is, the $\xi_j$ term of the $i$-th formula in \verb|X_formulas|. We obtain the coefficients and respective monomials of this $j$-th term, then for each (coefficient, monomial) pair given by \verb|coeff, monomial|, we extract the \verb|monomial_index| defined in \ref{collect_monomials}. We set the \verb|index_shift + monomial_index|'s entry of the sparse vector to be the monomial's coefficient; that is, we situate the monomial's coefficient within its respective $(\xi_i)$--box. Finally, we set the $i$-th column of the matrix to be this vector of coefficients we have constructed for each formula. 

\begin{exercise}
    Apply this construction to the example in this section, to see how the \verb|index_shift| is used to truncate the columns of the matrix for each $(\xi_i)$--box of monomial coefficients. 
\end{exercise}



\subsection{Collect linearly independent formulas and respective micro-graphs}\label{lin_ind}

We constructed the evaluation matrix in \S\ref{eval_matrix} which represents the coefficients of the monomials in the formulas for the micro-graphs. Each column of the matrix contains this information for one micro-graph. This means that by identifying those columns which span the column space of the evaluation matrix, we identify the micro-graphs whose formula spans the space of all other micro-graphs' formulas. So, we can establish linear independence of the formulas. Namely, we will use the command \verb|.pivots()| which returns the indices of the linearly independent columns of the graph-to-evaluation matrix. 

\begin{verbatim}
    pivots = X_evaluation_matrix.pivots()
    print("Maximal subset of linearly independent graphs:", list(pivots), 
            flush = True)
    lin_ind_graphs = [X_graphs[k] for k in pivots]
    lin_ind_flas = [X_formulas[k] for k in pivots]
\end{verbatim}

We now have ($i$) a list \verb|lin_ind_graphs| of the micro-graphs which produce linearly independent formulas, and ($ii$) the list \verb|lin_ind_flas| of linearly independent formulas themselves. Note that the indexing is preserved in both lists.

\begin{ex}
    Let us assume that we have a set of six micro-graphs from which we obtain six formulas, indexed $0,1,2,3,4,5$. Now let us assume that $0,2,4$ are linearly independent. These are the indices that will be returned in \verb|pivots|. Then we can call the suitable micro-graphs and formulas using \verb|X_graphs[k] for k in pivots| and \verb|X_formulas[k] for k in pivots|, respectively. 
\end{ex}\hfill$\blacksquare$ 

\begin{rem}
    It is not sufficient to take the $d$-dimensional graphs which produce linearly independent formulas in dimension $d$ and still obtain a trivialising vector field $\vec{X}_{d+1}^\gamma(P)$ in $(d+1)$, only searching over their descendants (see proposition 9 in \cite{IV}).
\end{rem}

\subsection{Skew-symmetrising the encodings (for dimension $d\geq4$) w.r.t. $a_i$}\label{skew}

We recall Lemma 5 from \cite{MSJB} concerning the property of $\vec{X}^\gamma_d(P)$ with regard to the swaps of Casimirs $a_1,a_2,...,a_{d-2}$. In 4D, this is:

\begin{lem} The Nambu--Poisson bracket $P(\varrho,a^1,a^2)$ is skew-symmetric under the swap $a^1\rightleftarrows a^2$: $$P(\varrho,a^1,a^2)=-P(\varrho,a^2,a^1).$$ The $\gamma_3$-flow $Q^{\gamma_3}_{d=4}(P)$ is built of four copies of $P$, therefore $Q^{\gamma_3}_{d=4}(P)$ is symmetric under $a^1\rightleftarrows a^2$; by swapping $a^1$ and $a^2$, we accumulate four minus signs: \begin{multline*}Q^{\gamma_3}_{d=4}\Big(P(\varrho,a^2,a^1)\otimes P(\varrho,a^2,a^1)\otimes P(\varrho,a^2,a^1)\otimes P(\varrho,a^2,a^1)\Big)\\=(-)^4Q^{\gamma_3}_{d=4}\Big(P(\varrho,a^1,a^2)\otimes P(\varrho,a^1,a^2)\otimes P(\varrho,a^1,a^2)\otimes P(\varrho,a^1,a^2)\Big)\\=Q^{\gamma_3}_{d=4}\Big(P(\varrho,a^1,a^2)\otimes P(\varrho,a^1,a^2)\otimes P(\varrho,a^1,a^2)\otimes P(\varrho,a^1,a^2)\Big).\end{multline*} Therefore, to find a vector field $\vec{X}^{\gamma_3}_{d=4}(P)$ such that $$\dot{P}=Q^{\gamma_3}_{d=4}(P)=\llbracket P,\vec{X}^{\gamma_3}_{d=4}(P)\rrbracket,$$ we need to find $\vec{X}^{\gamma_3}_{d=4}(P)$ which is skew-symmetric under $a^1\rightleftarrows a^2$. This can be seen by the fact that $\vec{X}^{\gamma_3}_{d=4}(P)$ is built of three copies of $P$, so accumulates three minus signs when swapping $a^1$ and $a^2$, therefore gives $(-)^3=-$, therefore is skew-symmetric under $a^1\rightleftarrows a^2$.
\end{lem}

To account for this skew-symmetry of $\vec{X}^{\gamma_3}_{d=4}(P)$ under the swap $a^1\rightleftarrows a^2$, we construct skew pair formulas $f$, and search for $\vec{X}^{\gamma_3}_{d=4}(P)$ over these: $$f=\tfrac{1}{2}\Bigl(\phi\bigl(\Gamma(a^1,a^2)\bigr)-\phi\bigl(\Gamma(a^2,a^1)\bigr)\Bigr).$$ To construct skew pairs, we take the formula of the graph $\Gamma$ with ordering of edges to Casimirs $a^1,a^2$ with $a^1\prec a^2$, and subtract the formula of the graph $\Gamma$ with ordering $a^2\prec a^1$. We divide by 2 to preserve the coefficients. By construction, each skew pair is purely obtained at the level of formulas. 

For sufficiently high $d\geq4$, the skew tuple formulas for finding $\vec{X}^{\gamma_3}_d(P)$ can be constructed by $$\text{ skew tuple }=\sum_{\sigma\in S_{\{1,2,...,d-2\}}\setminus\{id\}}(-)^{|\sigma|}\phi\bigl(\Gamma(\text{encoding}(\sigma(\boldsymbol{a})\bigr),$$ for $\boldsymbol{a}=\{a_1,a_2,...,a_{d-2}\}$. That is, for each micro-graph encoding $e$ which gives linearly independent formula $\phi(\Gamma(e))$, we produce new encodings $e(\sigma(a_i))$, running over almost all permutations $\sigma\in S_{\{1,2,...,d-2\}}\setminus\{id\}$; meaning we effectively permute the Casimirs $a_1,a_2,...,a_{d-2}$. 

\begin{ex}
    Take the 4D micro-graph encoding $$e=(0,1,4,7;1,3,5,8;1,2,6,9).$$ We have that $S_{\{1,2\}}\setminus\{id\}=\{\sigma=(12)\}$. In each permutation $\sigma$, object 1 corresponds to $a_1=4,5,6$, object 2 corresponds to $a_2=7,8,9$. Then $\text{encoding}\bigl(\sigma(a_i)\bigr)$ gives one new encoding: $$(0,1,7,4;1,3,8,5;1,2,9,6).$$
\end{ex}\hfill$\blacksquare$

Having skew-symmetrised a given set of micro-graph encodings, we repeat steps from \S\ref{encodings_to_mgs}--\S\ref{collect_monomials}. That is, we obtain the micro-graphs from the given set of skew-symmetrised encodings, we obtain the formulas of these micro-graphs, and we collect the differential monomials in these formulas.

\subsection{Repeat \ref{collect_monomials}--\ref{lin_ind} for skew pairs}\label{repeat}

Now, we want to repeat the following steps on the skew pair formulas: collect the monomials from the skew pair formulas (\ref{collect_monomials}), express the graph-to-formula evaluation as a matrix (\ref{eval_matrix}), and finally identify which of the skew pair formulas are linearly independent (\ref{lin_ind}). 

The code for these steps is the same as we had in \ref{collect_monomials}--\ref{lin_ind}. The objects we use are different. For example, instead of using \verb|X_formulas|, we use \verb|skew_pair_formulas|. The code for collecting the monomials from the skew pair formulas is:

\begin{verbatim}
X_monomial_basis = [set([]) for i in range(d)]
for i in range(d):
    for X in skew_pair_formulas:
        X_monomial_basis[i] |= set(X[i].monomials())
X_monomial_basis = [list(b) for b in X_monomial_basis]
X_monomial_index = [{m : k for k, m in enumerate(b)} for b in X_monomial_basis]

print("Number of monomials in components of skew_pair_formulas:", 
        [len(b) for b in X_monomial_basis], flush=True)
\end{verbatim}

\noindent The code for expressing the graph-to-formula evaluation as a matrix is:

\begin{verbatim}
X_monomial_count = sum(len(b) for b in X_monomial_basis)

X_evaluation_matrix = matrix(QQ, X_monomial_count, 
            len(skew_pair_formulas), sparse=True)
for i in range(len(skew_pair_formulas)):
    v = vector(QQ, X_monomial_count, sparse=True)
    index_shift = 0
    for j in range(d):
        f = skew_pair_formulas[i][j]
        for coeff, monomial in zip(f.coefficients(), f.monomials()):
            monomial_index = X_monomial_index[j][monomial]
            v[index_shift + monomial_index] = coeff
        index_shift += len(X_monomial_basis[j])
    X_evaluation_matrix.set_column(i, v)
\end{verbatim}

\noindent The code for identifying the linearly independent skew pair formulas is:

\begin{verbatim}
pivots = X_evaluation_matrix.pivots()
print("Maximal subset of linearly independent graphs:", list(pivots), 
            flush=True)
X_graphs_independent_skewed = [X_graphs[k] for k in pivots]
X_formulas_independent_skewed = [skew_pair_formulas[k] for k in pivots]
\end{verbatim}

Notice that in the list \verb|X_graphs_independent_skewed|, the indices identify graphs which were used to obtain the skew pair formulas. The graph themselves are neither skewed nor linearly independent as graphs. 

Now, \verb|X_formulas_independent_skewed| is the list of linearly independent skew pair formulas. We search for a solution to the $(\dot{a}_i,\dot{\varrho})$ system over these. 

\subsection{Set up $(\dot{a}_i,\dot{\varrho})$ system}\label{set up clever}

We know from \cite{skew21,skew23,MSJB} that if a vector field $\vec{X}$ satisfies the $(\dot{a}_i,\dot{\varrho})$ system 

\begin{subequations}
\begin{align}
\dot{\boldsymbol{a}} &=\llbracket \boldsymbol{a},\vec{X}\rrbracket\label{eqn:line-1} \\
\dot{\varrho}\cdot \partial_{\boldsymbol{x}} &= \llbracket \varrho\cdot \partial_{\boldsymbol{x}},\vec{X}\rrbracket \label{eqn:line-2}
\end{align}
\label{clever}
\end{subequations}

\noindent where \eqref{eqn:line-1} consists of $(d-2)$ equations and \eqref{eqn:line-2} consists of one equation; then $\vec{X}$ satisfies the coboundary equation $Q_d^\gamma=\llbracket P,\vec{X}\rrbracket$. To determine exactly the number of equations in the $\bigl( (d-2)+1\bigr)$ component algebraic systems, we track the number of differential monomials appearing when $\vec{X}$ acts on $\boldsymbol{a}$ $=(a_1,a_2,...,a_{d-2})$ or $\varrho$. Each linear algebraic equation obtained is a balance of the coefficient of one differential monomial in $\boldsymbol{a},\varrho$: namely, the coefficient in the known $\dot{a}_i$ or $\dot{\varrho}$ and the coefficient in $\llbracket a_i,\vec{X}\rrbracket$ or $\llbracket \varrho\partial_{\boldsymbol{x}},\vec{X}\rrbracket$. Using the $(\dot{a}_i,\dot{\varrho})$ system \eqref{clever} has two benefits compared with the coboundary equation $Q_d^\gamma=\llbracket P,\vec{X}\rrbracket$: it has less components, and the monomials are shorter. Namely, \eqref{clever} has $\bigl( (d-2)+1\bigr)$ sub-systems; the coboundary equation has $d(d-1)/2$ sub-systems; and monomials in $\llbracket P,\vec{X}\rrbracket$ are longer than those in $\llbracket a_i,\vec{X}\rrbracket$ or $\llbracket \varrho\partial_{\boldsymbol{x}},\vec{X}\rrbracket$.

Before we explicitly construct the $(\dot{a}_i,\dot{\varrho})$ system, we must import the implementation of the undirected graph complex. This is because the construction of $\dot{\boldsymbol{a}}$ and $\dot{\varrho}$ require the map $\Or(\gamma)$ and the bi-vector $Q^{\gamma}_d$, both of which are obtained via a cocycle $\gamma$ in the Kontsevich graph complex. Namely, $\gamma$ is a cocycle if it is an element of $\ker \mathsf{d}$, where $\textsf{d}=[\bullet\!\text{\textbf{--}}\!\bullet,\cdot]$ is the differential in the Kontsevich graph complex. 

\begin{verbatim}
    from gcaops.graph.undirected_graph_complex
    import UndirectedGraphComplex
\end{verbatim}

We take the undirected graph complex over the field of rationals $\mathbb{Q}$, with basis consisting of representatives of isomorphism classes of undirected graphs with no automorphisms that induce an odd permutation on edges; and store the graphs as collections of vectors, and the differentials as matrices. By differential, we mean $\textsf{d}=[\bullet\!\text{\textbf{--}}\!\bullet,\cdot]$ such that graphs on $V$ vertices and $E$ edges are mapped to a higher bi-grading, see \cite{JNMP2017}: $$(V,E)\xrightarrow[]{d}(V+1,E+1).$$ We choose to represent the differentials \textsf{d}, restricted to a given $(V,E)$--bi-grading, as above, as sparse matrices (meaning that most entries are zero), see tables 2 and 3 in \cite{JNMP2017} for details.

In this guide, we examine the coboundary equation $Q^{\gamma}_d(P)=\llbracket P,\vec{X}^{\gamma}_d(P)\rrbracket$ for the tetrahedral graph $\gamma_3$ in the graph complex. It is the 3-wheel cocycle built on 4 vertices and 6 edges. 

\begin{verbatim}
    GC = UndirectedGraphComplex(QQ, implementation = 'vector', sparse = True)
    tetrahedron = GC.cohomology_basis(4,6)[0]
\end{verbatim}

We have that \verb|GC.cohomology_basis(4,6)| spans the kernel of $\textsf{d}\vert_{(V=4,E=6)}$, which is a vector space naturally with a basis. The tetrahedron is the only admissible graph in $(4,6)$ and simultaneously a graph cocycle, making it the only element of a basis in $\ker \textsf{d}\vert_{(V=4,E=6)}$.

For other cocycles built on $V$ vertices and $E$ edges, the second line would become \begin{verbatim}
    name = GC.cohomology_basis(V,E)[i]
\end{verbatim} where $i$ can be equal to $0,1,...,b-1$, where $b$ is the number of elements in a basis of cocycles. We can choose and study our favourite (linear combination of) basis elements. 

We import the implementation of the directed graph complex.

\begin{verbatim}
    from gcaops.graph.directed_graph_complex
    import DirectedGraphComplex
\end{verbatim}

The program stores the graphs and differentials as collections of vectors, and the differentials as matrices, just like we had for the undirected graph complex. Note that \verb|DirectedGraphComplex| restricts the graphs to get rid of graphs with 1-cycles (loops), see p.132 of \cite{BuringPhD}. 

The canonical map $\Or$ (from \cite{ascona}, see \cite{OrMorphism2018}) from the undirected graph complex to the directed graph complex works as follows: every edge, independently from all others (if any), is directed consecutively in two opposite ways. This creates a sum of directed graphs. We implement this map as a conversion from GC to dGC, hence creating $2^6=64$ graphs in the linear combination \verb|tetrahedron_oriented|. 

\begin{verbatim}
dGC = DirectedGraphComplex(QQ, implementation = 'vector')
tetrahedron_oriented = dGC(tetrahedron)
tetrahedron_oriented_filtered = tetrahedron_oriented.filter(
                max_out_degree = 2)
tetrahedron_operation = Sd.graph_operation(tetrahedron_oriented_filtered)
\end{verbatim}

    The third line is to ensure that the directed tetrahedron is indeed built of wedges $\smash{\xleftarrow{L}\!\!\bullet\!\!\xrightarrow{R}}$ or single arrows $\smash{\bullet\!\!\xrightarrow{}}$ or vertices $\bullet$. This is because we know in advance that we will evaluate $\Or(\gamma_3)(\cdot,\cdot,\cdot,\cdot)$ at 4-tuples consisting of bi-vectors $P$ and scalar functions $a_i$, hence the multi-vectors at hand are at most bi-vectors. For instance, $Q^{\gamma_3}=\Or(\gamma_3)(P,P,P,P)$ is a bi-vector, and $\Or(\gamma_3)(P,P,P,a_i)$ is a 0-vector. The fourth line creates $\Or(\gamma_3)$. Specifically, $\Or(\gamma_3)(\cdot,\cdot,\cdot,\cdot)$ is an operation of arity 4 (we can plug in four multi-vectors) and degree $-6$ (it eats six $\xi_\alpha$ terms) on the Superfunction Algebra over the differential polynomial ring; the ring and the algebra were defined in \ref{graph_to_fla_eval}. The loss of six $\xi_\alpha$ can be seen by the fact that the six arrows each contain the following derivatives: $$(\text{tail})\frac{\overleftarrow{\partial}}{\partial\xi_i}\otimes\frac{\overrightarrow{\partial}}{\partial x^i}(\text{head}).$$ 

    We now have the operation $\Or(\gamma_3)(\cdot,\cdot,\cdot,\cdot)$, and can define the Casimir flows $\dot{\boldsymbol{a}}$: $$\dot{a}_i=4\cdot\Or(\gamma_3)(P,P,P,a_i),$$ where we let $P$ be the $d$-dimensional Nambu--Poisson bracket with pre-factor $\varrho$ and Casimir functions $a_1,a_2,...,a_{d-2}$. We implement the Nambu--Poisson bracket in terms of nested Schouten brackets $\llbracket\cdot,\cdot\rrbracket$, which can be expressed as follows, see \cite{skew21}: $$P=\llbracket \llbracket ... \llbracket \llbracket \varrho\xi_0\xi_1...\xi_{d-1},a_1\rrbracket,a_2\rrbracket,...\rrbracket,a_{d-2}\rrbracket,$$ and because \verb|.bracket()| refers to the Schouten bracket, we have:
    \begin{verbatim}
        P = (rho*epsilon).bracket(a1).bracket(a2)....bracket(a(d-2))
    \end{verbatim}
    We compute the Casimir flows $\dot{a}_i$: 

    \begin{verbatim}
def casimir_flow(f):
    return 4*tetrahedron_operation(P,P,P,f)

a = [Sd(a1), Sd(a2), ... , Sd(a(d-2))]

print("Calculating adot", flush=True)

adot = [casimir_flow(a_i) for a_i in a]

print("Calculated adot", flush=True)  
    \end{verbatim}
We compute $\llbracket X,a_i\rrbracket$ for each $X$ in \verb|X_formulas_independent_skewed| and each $i=1,2,...,d-2$. Notice the difference in ordering from \eqref{eqn:line-1}; this will lead to an extra minus sign at the end of the code. This can of course be adjusted for future implementations of the problem. 
\begin{verbatim}
print("Calculating X_a_formulas", flush=True)

X_a_formulas = [[X_formula.bracket(f) for X_formula in 
            X_formulas_independent_skewed] for f in a]

print("Calculated X_a_formulas", flush=True)
\end{verbatim}
Now, we can express a basis of the formulas $\llbracket X,a_i\rrbracket$, just like in \ref{collect_monomials}.
\begin{verbatim}
X_a_basis = [set(f[()].monomials()) for f in adot]
for k in range(len(a)):
    for X_a_formula in X_a_formulas[k]:
        X_a_basis[k] |= set(X_a_formula[()].monomials())
X_a_basis = [list(B) for B in X_a_basis]

print("Number of monomials in X_a_basis:", len(X_a_basis[0]), 
                    len(X_a_basis[1]), flush=True)
\end{verbatim}
We now define a graph-to-formula evaluation matrix as in \ref{eval_matrix}. The only difference with \ref{eval_matrix} is that we do not have to truncate the columns of the matrix to account for the $\xi_{i_1}$ terms. We express the linear map from graphs $g$ to formulas $\llbracket X_g,a_i\rrbracket$ given by $g\mapsto X_g(a_i)=\llbracket X_g,a_i\rrbracket$ as a matrix.
\begin{verbatim}
print("Calculating X_a_evaluation_matrix", flush=True)

X_a_monomial_index = [{m : k for k, m in enumerate(B)} for B in X_a_basis]
X_a_evaluation_matrix = [matrix(QQ, len(B), len(X_graphs_independent_skewed), 
                            sparse=True) for B in X_a_basis]
for i in range(len(a)):
    for j in range(len(X_graphs_independent_skewed)):
        v = vector(QQ, len(X_a_basis[i]), sparse=True)
        f = X_a_formulas[i][j][()]
        for coeff, monomial in zip(f.coefficients(), f.monomials()):
            monomial_index = X_a_monomial_index[i][monomial]
            v[monomial_index] = coeff
        X_a_evaluation_matrix[i].set_column(j, v)

print("Calculated X_a_evaluation_matrix", flush=True)
\end{verbatim}
We have computed the matrix representing the right hand side of \eqref{eqn:line-1}, $\llbracket \boldsymbol{a},\vec{X}\rrbracket$. We now express $\dot{\boldsymbol{a}}$ as a $(d-2)\times1$ vector.
\begin{verbatim}
print("Calculating adot_vector", flush=True)

adot_vector = [vector(QQ, len(B)) for B in X_a_basis]
for i in range(len(a)):
    f = adot[i][()]
    for coeff, monomial in zip(f.coefficients(), f.monomials()):
        monomial_index = X_a_monomial_index[i][monomial]
        adot_vector[i][monomial_index] = coeff

print("Calculated adot_vector", flush=True)
\end{verbatim}

We have implemented the ingredients necessary for the first line of the $(\dot{a}_i,\dot{\varrho})$ system \eqref{eqn:line-1}. We move on to implementing the ingredients necessary for the first line of the $(\dot{a}_i,\dot{\varrho})$ system \eqref{eqn:line-2}. To calculate $\dot{\varrho}$, we must first calculate the bi-vector $Q^{\gamma}_d(P)=\Or(P,P,P,P)$, because from \cite{skew21}, I., we have that $\dot{\varrho}$ is expressed from: \begin{equation}\label{rhodot}\dot{\varrho}\cdot\det\begin{vmatrix}\frac{\partial(f,g,a_1,...,a_{d-2})}{\partial(x^1,...,x^d)}\end{vmatrix}=\Bigl(Q^{\gamma}_d(P)-\sum_{i=1}^{d-2}P(\varrho,a_1,...,\dot{a_i},...,a_{d-2}\Bigr)(f,g).\end{equation} We construct the terms in the right hand side of \eqref{rhodot}, which we name \verb|Q_remainder|: 
\begin{verbatim}
print("Calculating Q_tetra", flush=True)

Q_tetra = tetrahedron_operation(P,P,P,P) 

#Note: this is a very costly operation, in time and memory! (To do: numbers??)

print("Calculated Q_tetra", flush=True)

print("Calculating rhodot", flush=True)

P0 = (rho*epsilon).bracket(adot[0]).bracket(a2).....bracket(a(d-2))
P1 = (rho*epsilon).bracket(a1).bracket(adot[1]).....bracket(a(d-2))
Pi=(rho*epsilon).bracket(a1).....bracket(adot[i]).....bracket(a(d-2))
P(d-3)= (rho*epsilon).bracket(a1).bracket(a2).....bracket(adot[d-3])
Q_remainder = Q_tetra - P0 - P1 - ... - Pi - ... - P(d-3)
\end{verbatim}
Now, we construct the terms on the left hand side of \eqref{rhodot}, which we name \verb|P_withoutrho|:
\begin{verbatim}
P_withoutrho = epsilon.bracket(a1).bracket(a2).....bracket(a(d-2))
\end{verbatim}
Finally, we can extract $\dot{\varrho}$. To do this, we notice that we only need to divide the $\xi_0\xi_1$ term on the right hand side by the $\xi_0\xi_1$ term on the left hand side. [To do: explain in more detail - how much detail?]
\begin{verbatim} 
rhodot = Q_remainder[0,1] // P_withoutrho[0,1]

print("Calculated rhodot", flush=True)
\end{verbatim}
We check that this step was executed correctly by verifying we indeed still obtain $Q^{\gamma}_d(P)$ as in \eqref{rhodot}:
\begin{verbatim}
print("Have nice expression for Q_tetra:", 
Q_tetra == rhodot*P_withoutrho + P0 + P1, flush=True)
\end{verbatim}
Now, we can calculate the formulas $\llbracket\varrho\xi_0...\xi_{d-1},X\rrbracket$ to obtain the right hand side of \eqref{eqn:line-2}. Like we did for $\llbracket X,a_i\rrbracket$, we construct a basis of monomials and compute the corresponding evaluation matrix.
\begin{verbatim}
print("Calculating X_rho_formulas", flush=True)

X_rho_formulas = [X_formula.bracket(rho*epsilon) for 
            X_formula in X_formulas_independent_skewed]

print("Calculated X_rho_formulas", flush=True)

X_rho_basis = set(rhodot.monomials())
for X_rho_formula in X_rho_formulas:
    X_rho_basis |= set(X_rho_formula[0,1,2,3,...,d-1].monomials())
X_rho_basis = list(X_rho_basis)

print("Number of monomials in X_rho_basis:", len(X_rho_basis), flush=True)

print("Calculating X_rho_evaluation_matrix", flush=True)

X_rho_monomial_index = {m : k for k, m in enumerate(X_rho_basis)}
X_rho_evaluation_matrix = matrix(QQ, len(X_rho_basis), 
            len(X_graphs_independent_skewed), sparse=True)
for j in range(len(X_graphs_independent_skewed)):
    f = X_rho_formulas[j][0,1,2,3,...,d-1]
    v = vector(QQ, len(X_rho_basis), sparse=True)
    for coeff, monomial in zip(f.coefficients(), f.monomials()):
        monomial_index = X_rho_monomial_index[monomial]
        v[monomial_index] = coeff
    X_rho_evaluation_matrix.set_column(j, v)

print("Calculated X_rho_evaluation_matrix", flush=True)
\end{verbatim}
Now we construct the left hand side of \eqref{eqn:line-2} as a vector $\dot{\varrho}$:
\begin{verbatim}
print("Calculating rhodot_vector", flush=True)

rhodot_vector = vector(QQ, len(X_rho_basis), sparse=True)
for coeff, monomial in zip(rhodot.coefficients(), rhodot.monomials()):
    monomial_index = X_rho_monomial_index[monomial]
    rhodot_vector[monomial_index] = coeff

print("Calculated rhodot_vector", flush=True)
\end{verbatim}

Finally, we now have all components of the $(\dot{a}_i,\dot{\varrho})$ system \eqref{clever} implemented in the code. We can solve the linear algebraic system. 

\subsection{Solve the $(\dot{a}_i,\dot{\varrho})$ system}\label{solve clever}

In \ref{set up clever}, we constructed the $(\dot{a}_i,\dot{\varrho})$ system as follows: $\dot{\boldsymbol{a}}$ and $\dot{\varrho}$ are stored as vectors, and $\llbracket X,a_i\rrbracket$ and $\llbracket \varrho\xi_0...,\xi_{d-1},X\rrbracket$ are saved as matrices. Now, we stack the vectors and matrices together to solve the equations in the $(\dot{a}_i,\dot{\varrho})$ system \eqref{eqn:line-1} and \eqref{eqn:line-2} simultaneously. We will use \verb|big_matrix| to refer to the right hand side of the $(\dot{a}_i,\dot{\varrho})$ system, and \verb|big_vector| to refer to the left hand side of the $(\dot{a}_i,\dot{\varrho})$ system. 

Note that on \textsf{SageMath}, we can immediately stack matrices on top of each other, but we choose to define the separate one-row matrices pertaining to each $\dot{a}_i$ and stack them for clarity. That is, we construct the big matrix of the $(\dot{a}_i,\dot{\varrho})$ system by stacking the rows of \verb|X_a_evaluation_matrix| with the one-row matrix \verb|X_rho_evaluation_matrix|. Recall from \ref{set up clever} that because of implementing $\llbracket X,a_i\rrbracket$ instead of $\llbracket a_i,X\rrbracket$, we obtain an extra minus sign. We account for this when defining the big vector. \textsf{SageMath} solves the system for us, and gives a solution vector \verb|X_solution_vector| (if it exists!) which consists of coefficients of linearly independent skew pair formulas. We can assemble the actual trivialising vector field \verb|X_solution| by multiplying the coefficients found in \verb|X_solution_vector| by their corresponding linearly independent skew pair formulas, and summing. 

\begin{verbatim}
print("Calculating X_solution_vector", flush=True)

big_matrix = X_a_evaluation_matrix[0].stack(X_a_evaluation_matrix[1])...
            ..stack(X_a_evaluation_matrix[d-3]).stack(X_rho_evaluation_matrix)
big_vector = vector(list(-adot_vector[0]) + list(-adot_vector[1]) + ... 
            + list(-adot_vector[d-3]) + list(-rhodot_vector))
X_solution_vector = big_matrix.solve_right(big_vector)

print("Calculated X_solution_vector which gives coefficients of 
            skew pair formulas", flush=True)
print("X_solution_vector =", X_solution_vector, flush=True)
print("Note that the coefficients given in the X_solution_vector are 
            from the linearly independent skew pairs, indeed len(X_solution_vector)=", len(X_solution_vector), flush=True)

X_solution = sum(c*f for c, f in zip(X_solution_vector, 
                    X_formulas_independent_skewed))

print("Verify that P.bracket(X_solution) == Q_tetra:", 
            P.bracket(X_solution) == Q_tetra, flush=True)
\end{verbatim}
In the last line above, we verify that the solution found to the $(\dot{a}_i,\dot{\varrho})$ system \eqref{clever} is indeed a solution to the coboundary equation $ Q^\gamma_d(P)=\llbracket P,\vec{X}^\gamma_d(P)\rrbracket$. Furthermore, we can find the number of parameters in the solution, and express the shifts of the solution as vectors relating to coefficients of linearly independent skew formulas.
\begin{verbatim}
print("Number of parameters in the solution:", big_matrix.right_nullity(),
                    flush=True)

print("Basis of kernel:", big_matrix.right_kernel().basis(), flush=True)
\end{verbatim}

With this, we conclude the code used to determine or verify if a given set of micro-graphs offer good formulas to construct a trivialising vector field $\vec{X}^{\gamma}_d(P)$ such that $Q^{\gamma}_d(P)$ is a coboundary, that is, $Q^{\gamma}_d=\llbracket P,\vec{X}^{\gamma}_d(P)\rrbracket$. 

\section{Finding Poisson 2-coboundaries}

The code which computes the solutions to the homogenized equation $\llbracket P,\Delta\vec{X}_d\rrbracket=0$, that is, the Poisson 2-coboundaries $\Delta\vec{X}_d$, is stored as ancillary files of \cite{fs}, and should be understandable after having read this paper. These Poisson 2-coboundaries are the shifts of any found solution $\vec{X}^{\gamma}_d(P)$, in the sense that $$\llbracket P,\vec{X}^\gamma_d(P)+\Delta\vec{X}_d\rrbracket=\llbracket P,\vec{X}^\gamma_d(P)\rrbracket+\llbracket P,\Delta\vec{X}_d\rrbracket =\llbracket P,\vec{X}^\gamma_d(P)\rrbracket,$$ meaning any found solution $\vec{X}^\gamma_d(P)$ is unique up to however many $\Delta\vec{X}_d$ exist. 


\section{Conclusion}

We restate the following optimisation point:

\begin{optimisation} The method to obtain formulas from graphs in \S\ref{graph_to_fla_eval} consumes large amounts of memory; this could be optimised. We believe that the following method would be efficient in reducing the amount of memory used. Compute the formulas one by one, and sum them together with undetermined coefficients instantly, so that similar terms are collected. We then have formulas of the form $(c_1+c_2)\cdot\text{term}$ instead of $c_1\cdot\text{term}+c_2\cdot\text{term}$. 
\end{optimisation}


\section{Acknowledgements}

The authors thank first and foremost Ricardo Buring for not just writing the software package \textsf{gcaops}, but for extremely helpful discussions and advice, as well as help with coding and optimising the code, for which they are deeply grateful. The authors thank the Center for Information Technology of the University of
Groningen for access to the High Performance Computing cluster, H\'abr\'ok.

\end{document}